\documentclass[5p,times]{elsarticle}
\usepackage[fleqn]{amsmath}
\usepackage{amssymb}
\usepackage{bbm}

\journal{}

\begin{document}
\begin{frontmatter}

\title{Personalized Stopping Rules in Bayesian Adaptive Mastery Assessment}

\author{Anni Sapountzi$^a$, Sandjai Bhulai$^b$, Ilja Cornelisz$^a$, and Chris van Klaveren$^a$}

\affiliation{organization={Vrije Universiteit Amsterdam},
			faculty={Faculty of Behavioral and Movement Sciences,},
            city={Amsterdam},
            country={The Netherlands}}

\affiliation{organization={Vrije Universiteit Amsterdam},
			faculty={Faculty of Sciences,},
			department = {Department of Mathematics,},
            city={Amsterdam},
            country={The Netherlands}}

\begin{abstract}
We propose a new model to assess the mastery level of a given skill efficiently. The model, called Bayesian Adaptive Mastery Assessment (BAMA), uses information on the accuracy and the response time of the answers given and infers the mastery at every step of the assessment. BAMA balances the length of the assessment and the certainty of the mastery inference by employing a Bayesian decision-theoretic framework adapted to each student. All these properties contribute to a novel approach in assessment models for intelligent learning systems. The purpose of this research is to explore the properties of  BAMA and evaluate its performance concerning the number of questions administered and the accuracy of the final mastery estimates across different students. We simulate student performances and establish that the model converges with low variance and high efficiency leading to shorter assessment duration for all students. Considering the experimental results, we expect our approach to avoid the issue of over-practicing and under-practicing and facilitate the development of Learning Analytics tools to support the tutors in the evaluation of learning effects and instructional decision making.
\end{abstract}


\end{frontmatter}

%\IEEEauthorblockA{$^1$Vrije Universiteit Amsterdam, Faculty of Behavioral and Movement Sciences \\
%$^2$Vrije Universiteit Amsterdam, Faculty of Science, Department of Mathematics \\
%Email: {annisapi@icloud.com, s.bhulai@vu.nl, i.cornelisz@vu.nl, c.p.b.j.van.klaveren@vu.nl, }}

\section{Introduction} \label{sec:introduction}
Assessment plays a vital role in learning. It provides both instructors and students with feedback on the overall effectiveness of their teaching or learning~\cite{assessment, lalley2009classroom}. Furthermore, it is considered a powerful tool for enhancing memory and building fluency~\citep{lalley2009classroom,learnermodelsstop,pelanek2017elo,glance2013pedagogical}. Computerized adaptive assessments are often operationalized in learning systems to measure the learning effects~\cite{peng2019personalized, martin2018computer,pardos2017big} (e.g., the outcomes after practicing). The effects are measured based on a sequential question administration process. Each time an answer is provided, a mathematical model (i.e., learner model) dynamically estimates the knowledge level of students~\cite{learnermodelsstop}. The information of that estimate is leveraged for a decision concerning an aspect of the assessment (e.g., the total number of items administered, skill difficulty) such that the model adapts to the individual's knowledge~\cite{pelanek2017elo, martin2018computer, kaser2016stop, elo2016applications}.

Much research has been done on mathematical models~\cite{learnermodelsstop, kaser2016stop, elo2016applications, fast, sapountzi2019dynamic, pelanek2018conceptual, kaserbayesian, dkt,bkt1993rules} and personalization~\cite{learnermodelsstop, glance2013pedagogical, peng2019personalized, pardos2017big, personalized, irtchen2005personalized, bktthreshold,dessi2019bridging}. Both are essential in online education as they attempt to provide resources tailored to the learner's needs while integrating interactions, skills and competencies of knowledge. The two main student modeling approaches include Knowledge Tracing (KT), which typically uses Bayesian models, and Item Response Theory (IRT), which typically uses latent factor models. The former is linked with developing a mastered performance and optimal instructional policies~\cite{kaser2016stop, anderson1993rules} for Intelligent Tutoring Systems. The latter has been studied mostly in efficient adaptive assessments and testing~\cite{elo2016applications,catweiss1984application,capklinkenberg2011computer} and the selection rules of optimal items~\cite{martin2018computer, irtchen2005personalized, catweiss1984application}. However, the distinction between adaptive testing and practicing models is not always clear. It is commonly based on whether the model accounts for knowledge level growth over the interaction with the system~\cite{learnermodelsstop, pardos2017big, pelanek2018conceptual}.

The question of whether someone has attained mastery is a fundamental question in learning sciences and Artificial Intelligence in Education~\citep{assessment, lalley2009classroom, bkt1993rules,masterydiagnosisoverview,binder2002fluency}. It lies at the heart of observation of performance and evaluation regarding the success of the precedent decisions~\citep{lalley2009classroom}. Assessment based on mastery is a category of adaptive assessment models designed to determine whether a student has attained a prespecified fluency level on a skill or subject ~\cite{lalley2009classroom, peng2019personalized, catweiss1984application, masterydiagnosisoverview, mastery1999bayesian, rt2019models, masteryrupp}. Mastery assessments are considered a first step toward what eventually might develop into fluent performance~\cite{lalley2009classroom,peng2019personalized}. 

An adaptive mastery assessment has multiple practical applications. Firstly, they can be leveraged as a post-practicing assessment to provide evidence of the practicing effects (e.g., the instructional policy improved the learning outcomes)~\cite{peng2019personalized, personalized}. Secondly, it can be utilized as a pre-screening tool prior to the student's practicing in the system to diagnose their attained knowledge so far on a topic.  Lastly, it could be used as a mid-practicing assessment to review skills composing a subject~\cite{lalley2009classroom} and test for prior knowledge from previous units or courses, instead of merely assuming it. The outcomes of the assessment support the stakeholders in the decision-making process regarding student progress in many ways (e.g., provide timely feedback to students, design remedial instructional policies).

The estimation of mastery both for adaptive testing~\cite{peng2019personalized, rt2019models, hshs2016use, hshs, 2002responsetime} and practicing~\cite{learnermodelsstop, pelanek2017elo, elo2016applications, capklinkenberg2011computer, rt2019models, gong2015towards} can be improved by utilizing the information from response times, the elapsed times between consecutive elements of a sequence. The speed of response has been studied extensively in educational psychology, typically with the dual-process theory~\cite{rt2019models, rtlearninglanuage, rtlearning}.  In the case of measuring mastered performance, a fast and accurate knowledge retrieval can better quantify student fluency under certain conditions~\cite{lalley2009classroom, peng2019personalized}. In offline learning, mastery is also considered to go beyond mere accuracy to include the pace, effort or speed of performance~\citep{lalley2009classroom,binder2002fluency, rtlearninglanuage, rtworkingmemorys}. However, the inclusion of response times for mastery diagnosis on a skill~\cite{masterydiagnosisoverview, rt2019models, masteryrupp} and adaptive practicing of declarative knowledge~\cite{learnermodelsstop, pelanek2017elo} has been largely absent in the literature.

Much less attention has been paid to the stopping criteria of mastery assessment models. Stopping criteria determine whether a student should continue practicing the skill, progress to a more advanced skill, or cover prior knowledge. A central problem in this line of research is selecting the stopping criterion that balances between the length of the assessment  (i.e., low measurement effort) and the risk of an inaccurate estimation of mastery~\cite{catweiss1984application, gong2015towards, lee2012impact, pelanek2017experimental, cat2019comparestop}. Within this framework, a certain enough mastery estimate should be achieved within a relatively short-length assessment avoiding any unnecessary over or under-practice~\cite{pelanek2018conceptual}. 

Stopping criteria constitute an integral component in an adaptive learning system~\cite{learnermodelsstop,kaser2016stop}. They influence the interpretation of the model results, the variance of the estimated mastery as well as the length of the process, and the follow-up policies if any~\cite{learnermodelsstop,pelanek2018conceptual,masteryrupp}. Recently, new evaluation methods and stopping policies have been proposed~\cite{learnermodelsstop,kaser2016stop,pelanek2018conceptual,rollinson2015predictive}, but the research is still limited. Stopping policies for mastery assessment models are even less studied~\cite{pelanek2018conceptual,masterydiagnosisoverview}.

Furthermore, research has shown that the assessment can terminate at a suboptimal time. That can be due to a population-based model~\cite{lalley2009classroom, pelanek2017elo, irtchen2005personalized, lee2012impact}, point-based estimates~\cite{mastery1999bayesian}, a decision evaluation with offline data~\cite{learnermodelsstop, rt2019models, lee2012impact, rollinson2015predictive}, and stopping rules that do not consider specific cases of student performances that have been observed in practice~\cite{kaser2016stop,pelanek2018conceptual,cat2019comparestop,rollinson2015predictive,catnewstop}. A model that solves for at least one of the above aspects is typically derived at the cost of data samples~\cite{elo2016applications, mastery1999bayesian} or model complexity and computation time~\cite{masterydiagnosisoverview, mastery1999bayesian, masteryrupp}.

In this paper, we are interested in inferring student mastery to support an efficient, criterion-referenced assessment. Our contribution is that we develop an individualized method to dynamically diagnose whether mastery is attained on a specific skill. For this purpose, we investigate the problem of stopping policies and mastery criteria for the item administration process. In our model, mastery is represented with a joint model of accuracy and response time on a single skill with homogeneous practicing items. The consequence of that model is an adaptive-length assessment illustrated in a decision-theoretic framework to account for certainty in the estimate. Specifically, a partially observed Markov decision process model (POMDP) with Bayesian updates has been employed to model the assessment. To the best of our knowledge, there is no other assessment model that satisfies simultaneously the following four properties:  (i) a model that infers an individual's mastery with a few observations, (ii) the exploitation of learner response times on practicing items, (iii) decision making in a dynamic manner, (iv) and adaptive length control tailored to the individual's performance, without comparing them with other student performances.

The paper is organized as follows. Section~\ref{sec:related_work} outlines the background and the basis of our contribution. Next, Section~\ref{sec:problem_formulation} presents our approach for adaptive mastery assessment, detailing the model components and how they are integrated into our method. Section~\ref{sec:numerical_experiments} illustrates the experimental evaluation we performed and describes the properties of the model. Finally,  Section~\ref{sec:conclusion} envisions the future directions and practical implications of our model, and it concludes the paper.

\section{Related work} \label{sec:related_work}
In this section, we draw different interrelated findings from the literature that relate to our research objective. Firstly, we frame what we consider as mastery in an assessment process and briefly describe models that incorporate the speed and the accuracy of the student response. Then we discuss related work on suboptimal criteria for the application of decision making in mastery assessment administration. Lastly, we explore mastery criteria in Intelligent Tutoring Systems (ITS) and Computerized Adaptive Tests (CAT). 

\subsection{Assessing mastery on a skill}
Mastery has several related definitions. It has been often viewed as the phenomenon of fluency that succeeds ability \cite{peng2019personalized, binder2002fluency,rtlearninglanuage}. From that perspective, mastery is considered the third phase of the learning process~\cite{NOericsson2006influence}, where knowledge is automated and well-composed~\cite{lalley2009classroom, learnermodelsstop, pelanek2017elo}. Mastery is then viewed as ease of information processing. The individual has developed and stored the knowledge component required for the skill in working memory; and they can retrieve it from working memory with minimal effort and execute it accurately~\citep{learnermodelsstop, binder2002fluency, rtlearninglanuage,rtlearning,rtworkingmemorys}. 

A skill is typically assimilated to a knowledge component in learner models. It is tightly defined in Knowledge Tracing in contrast to latent factor models that model knowledge on a broad skill.  In the former case, a skill is composed of many independent, homogeneous problem-solving items. Then, the probability of demonstrating proficient performance to an item that involves the skill depends on the student mastery of the skill. Conversely, skill mastery is measured by the students' ability to recall the skill and apply it to solve any item that involves the skill.

\subsection{Response time in performance modeling}
A student who demonstrates mastery can be modeled as a student who never responds incorrectly~\cite{kaser2016stop}, but that would not be adequate for a mastery assessment~\cite{lalley2009classroom, learnermodelsstop, peng2019personalized, binder2002fluency, pelanek2017experimental}. 

For illustration, consider two specific contrasting examples, where the skill is the multiplication of single-digit numbers, and an item that assesses the solution of `$6 \times 4 = \,?$' . Students have practiced and are able to provide the correct answer. Further on their practicing, the student's knowledge reaches a point that without any further thoughts solves it to 24, rather than having to compute $6 \times 4$ is $4+4+4+4+4+4$. Given that a student will provide a correct response with high probability in both cases, a model based only on the information of correctness of the answer will not distinguish between these differential levels of mastery.

Response times have been used in modeling student performance typically with the objective of strengthening the estimation of attained mastery~\cite{learnermodelsstop, pelanek2017elo, peng2019personalized, rt2019models, hshs2016use, 2002responsetime}. They constitute a proxy for the amount of effort a person puts into a task~\cite{pelanek2017elo, elo2016applications, capklinkenberg2011computer, 2002responsetime, gong2015towards, rt2011analysis}, and they are considered a significant source of individual differences~\cite{pelanek2017elo,elo2016applications, capklinkenberg2011computer, rt2019models, hshs2016use, 2002responsetime, gong2015towards, hofman2018fast}. 

In the psychology of learning, response times are typically interpreted using the dual-process theories of thinking~\cite{2002responsetime,rtlearning} distinguishing between implicit knowledge -- involves automatic processing such as knowledge retrieval -- and explicit knowledge -- entails controlled processing such as sequential reasoning steps~\cite{rtlearninglanuage}. In the first case, the problem-solving items assess declarative, factual, fined-grained knowledge that does not require complex or procedural knowledge, similar to the example we provided above. Differently said, the person solving an item can usually be assumed not to be performing any other problem-solving tasks simultaneously. The items also typically have a closed-ended form (i.e., one possible correct answer to an item that is precisely given to the system). Assuming local dependency of response time and response accuracy on an item, the performance on the item can then be an indirect indication of the main type of processing: automated vs. controlled processing~\cite{rt2019models,rtlearninglanuage}. 

A well-known model that relates response time to accuracy is the speed-accuracy tradeoff~\citep{capklinkenberg2011computer, rt2019models, hshs2016use, hshs}. It has been previously used in educational systems, such as adaptive practicing of factual knowledge and assessments~\cite{learnermodelsstop, pelanek2017elo, capklinkenberg2011computer, hshs, pelanek2017experimental}. We use a similar rule to that, in the sense that we also favor a quick response with higher mastery. 

\subsection{Suboptimal stopping criteria}
In general, the rule that determines the termination of a dynamic-length process is called a stopping rule. A stopping policy refers to a sequence of time-dependent stopping rules. Student mastery profiles are distinguished based on a mastery criterion, a stopping rule that relates to a process involved in reaching or not a pre-defined mastery level~\cite{pelanek2018conceptual, catweiss1984application, masterydiagnosisoverview, rt2019models, masteryrupp}. The pre-defined value of the criterion ensures that all individuals are assessed to the same established quality level or instructional goal~\cite{lalley2009classroom, learnermodelsstop, binder2002fluency, cat2019comparestop}. 

Regardless of the application, an optimal criterion is offered by maximizing the probability of making a correct decision while simultaneously minimizing the length of the process~\cite{kaser2016stop, mastery1999bayesian}.  Specifically, in learning systems it has been reported problem-related properties that can lead to suboptimal decisions.  Firstly, a criterion is set on a fixed-length assessment is a disadvantage on a precise estimate across a range of different mastery levels. This is in contrast with a variable-length assessment that typically achieves a certain degree of precision for everyone, even at the cost of providing more items to some students compared to others~\cite{kaser2016stop, cat2019comparestop, catnewstop, stopcatleroux2014comparison}. 

Secondly, a criterion based on mastery estimates which do not explicitly quantify uncertainty has the potential to result in suboptimal stopping~\cite{pelanek2018conceptual}. In reality, there is a lot of uncertainty involved since performance is not a perfect information indicator of mastery (typos, guessing, or slipping behaviors~\cite{rt2019models}), latent sequential effects between a student's performance~\cite{kaser2016stop}, differences on prior student length of practice on the system (i.e., attrition bias)~\cite{learnermodelsstop}, or student mental state while taking the assessment. Models of Bayesian Knowledge Tracing involve a framework known in AI literature as decision making under uncertainty. That framework naturally retains the whole distribution of the process and maintains all available information on those estimates, thereby allowing it to stop according to the associated uncertainty. 

Furthermore, theoretically, a policy evaluation that uses online data (as opposed to offline) and a dynamic model (e.g., online estimates) is a more proper methodology, when the goal is to estimate whether the process should stop~\cite{learnermodelsstop}. That is to continuously update the estimates after observing each answer to prevent the introduction of specific biases~\cite{learnermodelsstop,pelanek2017elo}. On the other side, online estimations of a population-based model use the same data to re-train the model, which is prone to a feedback loop problem. In any case, such a methodology entails difficulties in real-world settings. 

In assessment theory, an individualized model is a necessary component -- assess each individual performance independent of other student performances -- and can be used for formative assessment and timely feedback~\cite{lalley2009classroom,binder2002fluency}. Hence, a population-based model includes a fundamental issue. It assumes that the initial probability for mastery or the perception of the `easiness' of a skill does not depend on the individual~\citep{kaser2016stop, rt2019models}. In practice, this has been shown to provide suboptimal results~\citep{lee2012impact, categgen2006optimal}. 

A common problem in intelligent tutors and adaptive computerized testing is that the policy may never stop for a specific class of individuals (e.g., wheel-spinners, low-engaged), and it systematically results in the over-practice problem. That occurs when further item administration will not provide a more certain or a higher point-based mastery estimate~\cite{learnermodelsstop, kaser2016stop,pelanek2018conceptual, catnewstop}. 

Additionally, the inclusion of student response times has been suggested to have a higher impact on mastery estimates compared to using a more complex model that additionally estimates item properties~\cite{learnermodelsstop,pelanek2017experimental} on decisions regarding mastery.

\subsection{Criteria in adaptive practicing and testing}
The motivation for stopping policies has become strongly present in learner models for adaptive practicing~\cite{kaser2016stop,pelanek2018conceptual, lee2012impact, pelanek2017experimental}. As the modeling choices on criteria have significant consequences on the length of the practicing process and the students' final mastery interpretation~\cite{learnermodelsstop, masterydiagnosisoverview}.

A conceptual tool for supporting the choice and interpretation of  mastery criteria has been recently proposed~\cite{pelanek2018conceptual}. The criterion is typically interpreted as a weighted proportion of correct answers over the total number of items~\cite{learnermodelsstop, pelanek2018conceptual, catweiss1984application, mastery1999bayesian} (latent factor models), or the uncertainty of the estimate for a fine-grained skill (Bayesian models)~\cite{pelanek2018conceptual, anderson1993rules}, or being modeled with a simple rule of a fixed number of consecutive correct answers~\cite{gong2015towards, pelanek2017experimental} (simple models on exponential moving average). For the latent factor models, the threshold is set to a 50\%--70\% probability answering correctly, depending on the level of the assessment's leniency. For the Bayesian models, the threshold is set to a value of 95\% certainty; and the latter ranges from 3--5 responses. 

More recent work introduces model-independent evaluation metrics \cite{teal} and policies termed as `when to stop' policies~\cite{kaser2016stop,rollinson2015predictive} compatible with several adaptive practicing models from the literature. These works are focused on addressing the `mastery attrition bias' which is caused by the non-random differences in the students' practicing length)~\cite{learnermodelsstop}  and evaluating the impact of a learner model regarding the length of the session and the accuracy of the estimate. The Teal metric defines a threshold for the balance between two processes: the ratio of correctly solved items -- after the policy stops -- and the number of total items -- before the policy stops. Generally, `when to stop' policies~\cite{rollinson2015predictive} propose to stop when the predicted probability that the student will respond correctly to the next item is not changing~\cite{rollinson2015predictive}:
\begin{equation}
\mathbbm{P}( | \, \mathbbm{P}(C_t) - \mathbbm{P}(C_{t+1}) \, | < \epsilon) >\delta,
\end{equation}
where $\mathbbm{P}(C_t)$ denotes the probability of observing a correct response at time $t$, $\delta$ is the fixed change set to $95\%$, and $\epsilon$ is set to a small value $0.01$. For the evaluation of the model, the maximum number of items was set to 25. An extension of that policy is recently proposed to consider wheel-spinning students~\cite{kaser2016stop}. However, this research does not focus on the response times of students, instead, the effort is modeled with the number of total answered items. 

IRT approaches have also investigated the stopping criteria in order to balance the dual concerns of measurement precision and testing efficiency~\cite{catweiss1984application, mastery1999bayesian, cat2019comparestop, catnewstop, stopcatleroux2014comparison}. However, they translate mastery regarding item properties assuming one broad skill with mastery being modeled as a continuous variable, similarly to latent models in adaptive practicing. To prevent the needless administration of items to wheel-spinning students, combined with a minimum or maximum number of items rule~\cite{kaser2016stop, cat2019comparestop, catnewstop}. 

Cognitive diagnosis models are classification models for a binary decision of mastery. They have gained attention for modeling fine-grained skills, but they are not yet present in online learning systems because of their high cost of requirements~\cite{masterydiagnosisoverview, masteryrupp}. In mastery assessments that use item response theory and cognitive diagnosis models, the threshold is usually set around $0.5$~\cite{catweiss1984application, mastery1999bayesian, masteryrupp}. 

To the best of our knowledge, no model regardless of whether they measure certainty or they are robust to data noise (e.g., wheel-spinning behaviors), accounts both for the information of response times and individualized estimations. In addition, we desire to design a system that (i) takes into account the properties that have been shown to lead to suboptimal stopping policies, (ii) assesses each student individually, where the initial mastery of a skill and the prior practice length to be dependent on an individual,  (iii) chooses a model that does not require many data samples, and it is suitable for online applications, (iv) measures certainty explicitly, (v) and uses a sequence of online estimates. We intend to bridge these gaps in one model. To address that, we frame the problem of stopping administering items within the artificial intelligence field of decision making under uncertainty.

\section{Problem formulation} \label{sec:problem_formulation}
In this section, we introduce and formalize the mastery assessment process. We focus only on the case of learning for a single knowledge component.

In our setting, there is an agent who assesses a single student. The student is presented items (e.g., exercises or assignments) in sequence. Each item has to be solved within a certain time. After the student has provided an answer, the agent observes if the answer is correct, and the time it took to provide the answer. After this observation, the agent has to decide whether the assessment should stop because the student has achieved mastery for this knowledge component, or whether to continue.

The agent is faced with the problem of how to convert the observations that he receives to statements on the mastery assessment. A correctly given answer contributes more to the mastery if it is given rapidly in a consistent manner. Hence, the correctness and the answering speed should be looked at together, rather than independently. Moreover, answering in a consistent manner also means that averages do not provide the complete picture, but one also needs to take into account the variance, or even other metrics.

Let $d$ denote the maximum response time that is allowed in answering an item. We denote the accuracy by $P \in \{0, 1\}$, where 0 indicates an incorrect answer, and 1 a correct answer. Furthermore, let $T \in \mathbbm{R}_+$ denote the response time of the student. We combine these two parts into a single score $Z$, similar to previous literature~\citep{capklinkenberg2011computer}, by
\begin{equation} \label{eq:metric}
Z = P \cdot \left(1- \frac{T}{d} \right)^+.
\end{equation}
The score $Z$ is close to 1 when a student answers correctly and relatively fast with respect to $d$. As the student delays the response, the value of $Z$ decreases. The score value becomes zero in two cases: when a student answers incorrectly, or when the response time exceeds $d$. 

The expression in Equation~(\ref{eq:metric}) implies that $Z$ is a deterministic quantity. However, in practice, the student will answer a fraction $\theta$ of the items correctly with response time $1/\lambda$ on average. In our model, this is reflected by saying that $P$ can be seen as a Bernoulli trial with parameter $\theta$, and that $T$ has an exponential distribution with rate $\lambda$. Using the Bernoulli distribution is common practice in modeling responses on a closed-form type of similar items~\citep{masterydiagnosisoverview, mastery1999bayesian}. The exponential distribution has been used in the past to model response times for factual knowledge, such as mental rotation items assuming sequential processes among items solving~\citep{rt2019models}. The instantaneous knowledge retrieval can be then considered as almost memoryless. Therefore, a student can be characterized by the pair  $(\theta, \lambda)$ for their performance on this single knowledge component. If the agent knows $(\theta, \lambda)$ with certainty, we can assess the student directly. Yet, in practice, this information is not known, and needs to be estimated from the observations. 

In this paper, we adopt a Bayesian approach to estimating the true unknown parameters $\theta$ and $\lambda$ of a student. For this purpose, we model $\theta$ and $\lambda$ itself as a random variable. This makes $P$ and $T$ doubly stochastic; namely, it is a random variable that has as parameter a random variable. To keep the formulation tractable, we model $\theta$ by a Beta distribution, and $\lambda$ by a Gamma distribution. The Beta distribution is the conjugate family of distributions for a Bernoulli trial, as is the Gamma distribution for exponential response times. The advantage of using conjugate families of distributions is that the posterior updates of $\theta$ and $\lambda$ will remain a Beta and Gamma distribution, respectively. Hence, in our problem, it suffices to only keep track of parameters of the distributions, instead of keeping track of a complete distribution that is difficult to characterize.

We suppose that our prior Beta distribution has parameters $\alpha$ and $\beta$. Similarly, we suppose that our prior Gamma distribution has parameters $n$ and $\gamma$. Now assume that a student has answered an item. Let $p$ denote the correctness of the response, which is either incorrect ($p = 0$) or correct ($p = 1$), and is given after $t$ time units. Then the posterior distributions are given by a 
\begin{equation} \label{eq:update_beta}
\text{Beta}(\alpha + p, \beta + 1 - p)
\end{equation}
distribution, and a 
\begin{equation} \label{eq:update_gamma}
\text{Gamma}(n+1, \gamma + t)
\end{equation}
distribution, respectively. As more responses become available, the posterior distributions of the accuracy (the Beta distribution) and the response time (the Gamma distribution) become more centered and peaked around the true values of $\theta$ and $\lambda$. Moreover, the agent can now utilize all information of the student, not only the averages, since the agent knows the complete distribution.

The probability of a correct response can be calculated as follows for a Beta($\alpha$, $\beta$) prior.
\begin{equation} \label{eq:beta}
\begin{split}
\mathbbm{P}(P=1) & = \int_0^1 \theta \cdot f_\theta (\theta;\alpha, \beta) \, \text{d}\theta 
\cr & = \int_0^1\frac{\Gamma(\alpha + \beta)}{\Gamma(\alpha) \cdot \Gamma(\beta)} \cdot \theta^\alpha \cdot (1-\theta)^{\beta-1} \, \text{d}\theta = \frac{\alpha}{\alpha + \beta},
\end{split}
\end{equation}
where $\Gamma$ is the Gamma function. Consequently, it follows that  $\mathbbm{P}(P=0)= \frac{\beta}{\alpha + \beta}$. 
Similarly, the response time density function evaluated at $t$ time units, for a Gamma($n$, $\gamma$), is given by:
\begin{equation} \label{eq:gamma}
\begin{split}
 f_T (t) & = \int_0^\infty \lambda e^{-\lambda t} \cdot f_\lambda (\lambda;n, \gamma) \, \text{d}\lambda = \int_0^\infty \frac{\gamma ^n \lambda ^n e^{-\lambda(\gamma + t)}}{\Gamma(n)} \ \text{d}\lambda 
 \cr & = \frac{n}{\gamma + t} \left(\frac{\gamma}{\gamma +t} \right)^n.
\end{split}
\end{equation}
Note that Equation~(\ref{eq:update_beta}) shows that the parameters $\alpha$ and $\beta$ could be given the interpretation of the number of correctly answered items and the number of incorrectly answered items, respectively, if the prior could be initialized with $\alpha = \beta = 0$. Thus, $\theta$ is estimated by the fraction of correctly answered items over the total number of administered items, see Equation~(\ref{eq:beta}). Of course, $\alpha$ and $\beta$ should be initialized differently (and need not be integer-valued) based on past information that is available on the students. Similarly, Equation~(\ref{eq:update_gamma}) shows that $n$ and $\gamma$ count the number of items and the total sum of response times, respectively, when the prior could be initialized with $n = \gamma = 0$. Therefore, $1 / \lambda$ is estimated by $\gamma / n$, which is the average response time observed from data, see Equation~(\ref{eq:gamma}).

\subsection{Personalized stopping rules}
Now that we have defined the ingredients of the model, i.e., the way the parameters $\theta$ and $\lambda$ are dynamically estimated, we can proceed to derive optimal stopping policies. The idea is to use the ingredients to adaptively administer the number of items to each student according to their estimated mastery level. This policy can be derived efficiently within the framework of partially observed Markov Decision Process (POMDP) models. To formulate the problem into a POMDP framework, one needs to define states $s$, admissible actions $a$, a reward function $r$, and a transition probability function $p$ for the assessment problem.

We define an information state $s = (\alpha, \beta, n, \gamma)$. This state represents the parameters of the Beta distribution and the Gamma distribution. This essentially serves as a proxy for information on $\theta$ and $\lambda$ that cannot be observed directly. Hence, the state serves as a sufficient statistic of the past answers over the assessment. 
Here, conjugated priors simplify the state description; instead of storing the complete distribution, one can suffice with keeping track of the parameters of the distribution only.

The actions $a$ are modeled by binary values; $a = 0$ (stop the assessment) and $a = 1$ (continue the assessment). The action is based on the current information of the system, i.e., the state $s$. Based on the state, we want to have an algorithm that can estimate the value of the $Z$ score accurately (e.g., the mean value), while at the same time take into account the uncertainty (e.g., the variance). Once the state $s$ has been observed, and an action $a$ has been chosen, the evolution of state $s$ to state $s'$ can be calculated using Equations~(\ref{eq:update_beta}) and~(\ref{eq:update_gamma}). Therefore, the transition function $p(s, a, s')$ is deterministic according to these equations.

Next to the evolution of state $s$ to state $s'$, the system also gives a signal in the form of a reward function $r$. The reward function $r$ represents a student's $Z$ score, as estimated by the distributions of $P$ and $T$. Thus, e.g., Equation~(\ref{eq:beta}) represents the probability that a correct answer is given, and Equation~(\ref{eq:gamma}) represents the density that an answer is given at $t$ time units. Therefore, $r(s, a, s') = \left(1- \frac{t}{d} \right)^+$ with the corresponding probabilities.

The expected value of the reward is discounted with a hyperparameter $\eta \in [0, 1]$. This describes how the rewards are weighted over time. For a value of $\eta$ close to 1, the algorithm puts more emphasis on exploration implying a potentially longer assessment, but with higher certainty and less variance. For a value of $\eta$ close to 0, the algorithm puts more emphasis on exploitation implying higher efficiency and more variance, i.e., the algorithm concludes quickly after a few offered items whether the student has passed the mastery level or not. In practice, the value of $\eta$ is chosen between 0.7 and 0.8 for a good balance between efficiency and certainty.

Now that all ingredients have been defined, the decision of stopping the assessment can be calculated recursively using the Bellman equations:
\begin{equation} \label{eq:bellman}
V(s) = \max_a \sum_{s'}  p(s, a, s') [r(s, a, s') + \eta V(s')],
\end{equation}
where the value function ${V(s)}$ measures the quality of state $s$. Given a predetermined parameter $\xi$, that describes the targeted mastery level, it can be determined in each state $s$ whether another item should be offered to the student. By substitution of Equations~(\ref{eq:metric})--(\ref{eq:gamma}) into Equation~(\ref{eq:bellman}), we obtain
\begin{equation} \label{eq:vi}
\begin{split}
V(\alpha, \beta, n, \gamma) = & \min \Bigg\{\int_0^\infty \Bigg(\frac{\alpha}{\alpha + \beta} \Big[\Big(1 - \frac{t}{d} \Big)^+ \ + 
\cr & \eta V(\alpha+1, \beta, n+1, \gamma+t) \Big] \ +
\cr & \frac{\beta}{\alpha + \beta}\eta V(\alpha, \beta+1, n+1, \gamma+t)\Bigg) \  \cdot
\cr & \frac{n}{\gamma+t}\Big(\frac{\gamma}{\gamma +t}\Big)^{n} \text{d}t \ ; \ \frac{\xi}{1-\eta} \Bigg\}.
\end{split}
\end{equation}
Note that there are a few differences with respect to Equation~(\ref{eq:bellman}). The variable $\gamma$ is a continuous variable, therefore, the sum in Equation~(\ref{eq:bellman}) is replaced by an integral. Moreover, the maximization in Equation~(\ref{eq:bellman}) is replaced with a minimization so that the algorithm stops the assessment once the student has obtained mastery. To illustrate that, imagine having two students with different performances. One student is being tested, and the other has the performance fixed at the threshold level $\xi$. When the algorithm chooses the second student, then this means that the first student has surpassed the threshold $\xi$ (since the minimum is chosen). 

\begin{figure}[!tb]
    \centering
    \includegraphics[width = \linewidth]{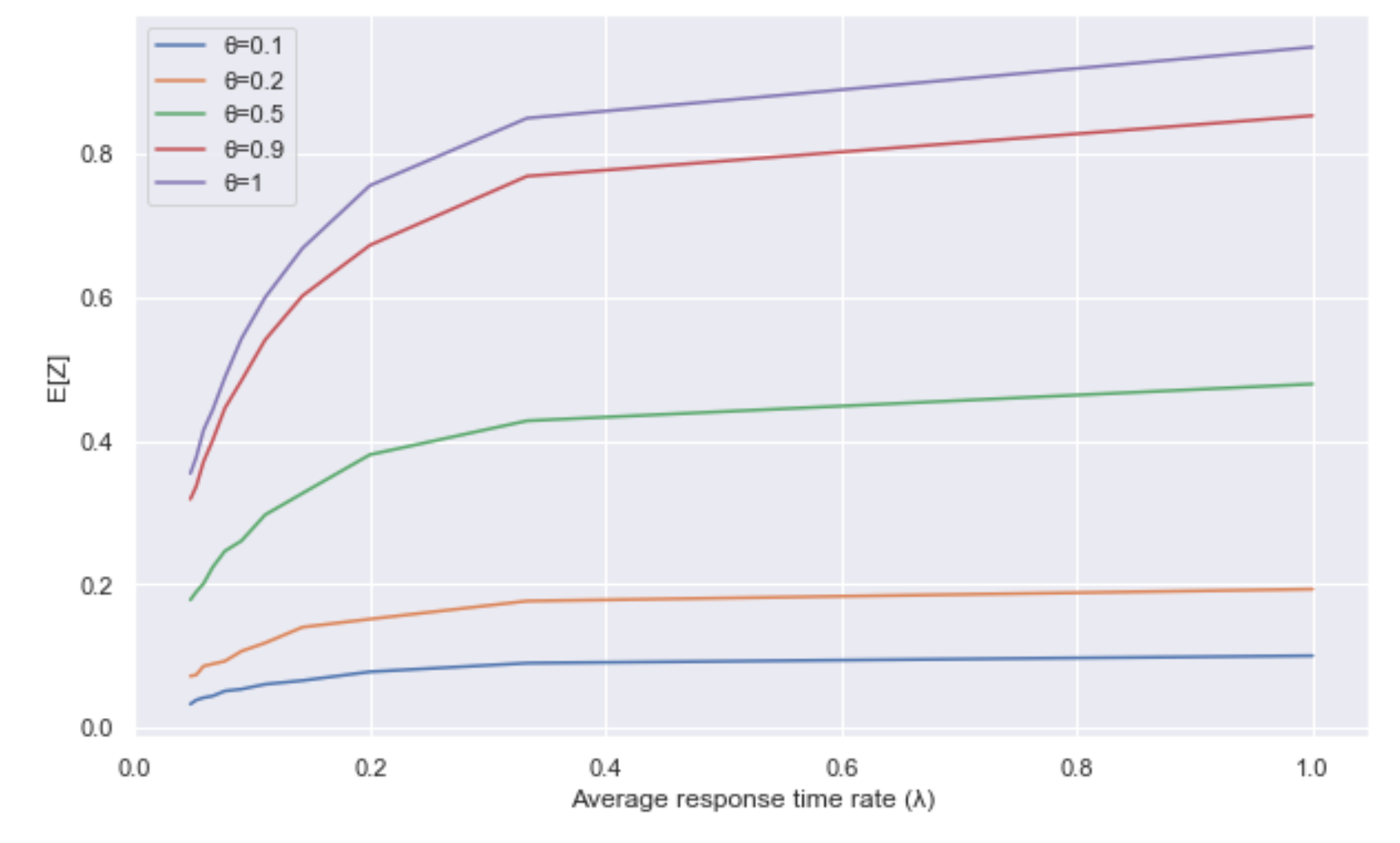}
    \caption{The relationship between the response time rate $\lambda$ and $\mathbbm{E}[Z]$. The quickest responses are given for $\lambda = 1$. The upper (purple) and lower (blue) lines show the highest and lowest score, respectively.}
    \label{fig:1}
\end{figure}

A skill is considered to be mastered in state $s$, once the first term in the minimization exceeds $\xi /(1 - \eta) $. Then, the agent stops the assessment with a number of total administered items that has been adapted to the student's mastery level. If the first term in the minimization does not exceed the latter, $\xi/(1 - \eta)$, then it is considered that the student has not mastered the skill. In that case, the student may need to solve all the items composing the skill. Alternatively, a simple decision rule based on existing literature can be used to stop, such as the 3 consecutive answers in a row of which not all are correct~\cite{pelanek2017experimental}. The literature regarding the stopping criteria of an assessment, as discussed in the previous section, together with the results of the experiments, will be used for the choice of the value for  $\xi$. Similar to all assessment models, the true mastery level of a student is never explicitly derived from the sequence of answers due to the unavailability of the ground truth (i.e., mastery is a hypothetical construct). 

\section{Numerical experiments} \label{sec:numerical_experiments}
In this section, we present three general numerical experiments to illustrate our model. These include the properties of the $Z$ score (i.e., model behavior), the convergence properties of the model (i.e., the evolution of the parameters over time), and the tuning of the predefined threshold as a mastery criterion. 

It is important to distinguish the parameters contained in the first experiment from the second experiment. The former simulates a student profile for a single skill as described by the combination of $(\theta,\lambda)$ values. The goal of the second experiment is to simulate students as they are consecutively responding to items within the assessment process. This is described by the parameters $\alpha, \beta, n, \gamma$. In the first experiment, we consider a skill with, in theory, an infinite number of items for each student, while in the latter, we constrain the set to 30 items.  For all experiments, we set the value of the maximum permitted response time to $d =20$, and the fastest answer to be $\lambda = 1$.

\subsection{Effect of accuracy and speed on mastery}
This section explores the sensitivity of the distribution $Z$ with regard to its parameters $(\theta,\lambda)$ such that the effect of the accuracy and the response time on the score $Z$ is studied. In practice, this illustrates how $Z$ behaves, given varying mastery profiles, at the time either before or after the assessment has taken place. In all figures, the concentration of $Z$ values around $0$ is due to students having an inaccurate or a slower than the permitted time response. As expected, such a pattern is penalized by the model with a score close to 0. 

\begin{figure}[!tb]
    \centering
    \includegraphics[width = \linewidth]{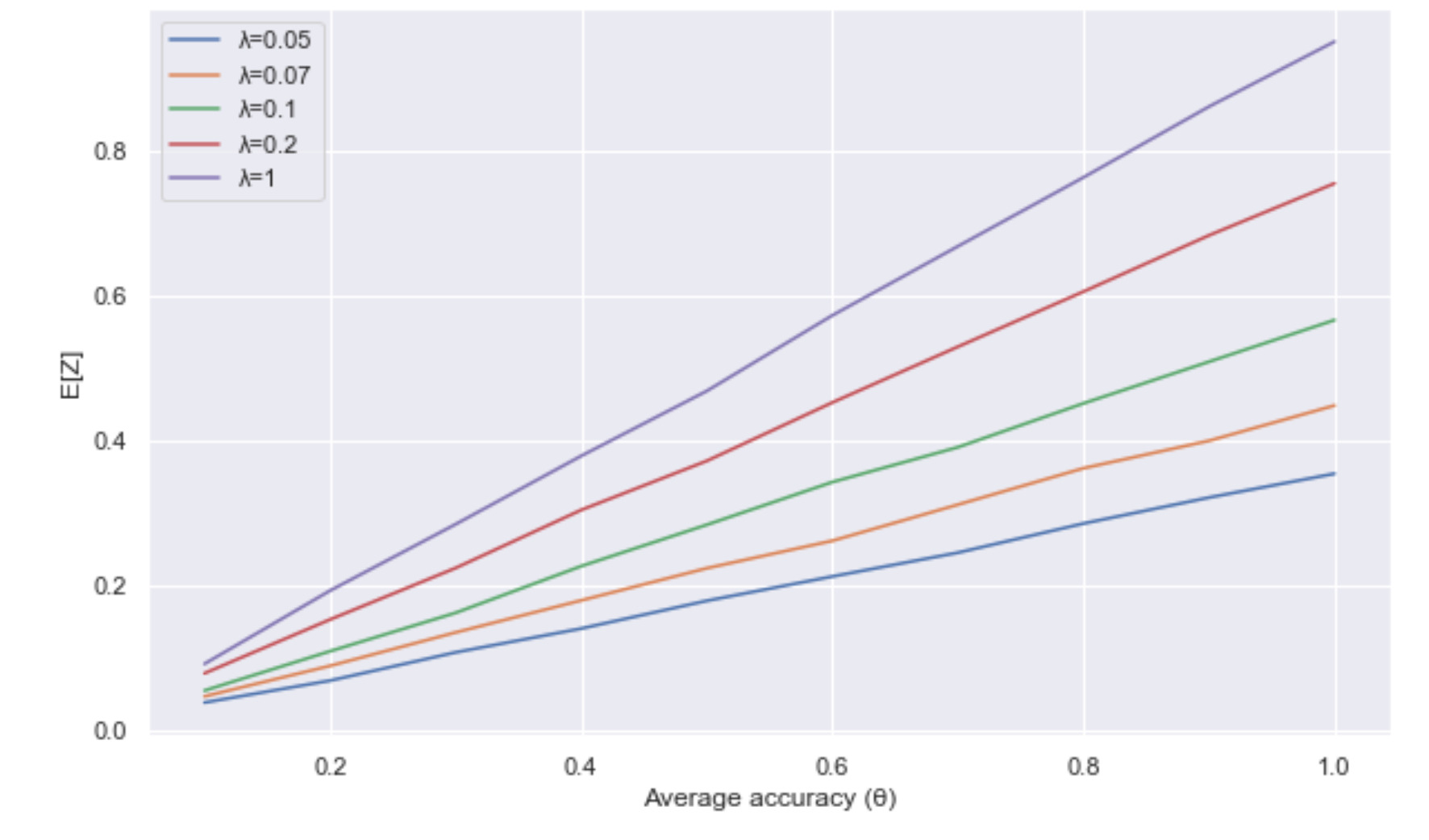}
    \caption{The relationship between the accuracy parameter $\theta$ and $E[Z]$ is linear with an upward slope. All responses are correct for $\theta = 1$. The upper (purple) and lower (blue) lines show the highest and lowest score, respectively.}
    \label{fig:2}
\end{figure}

Figures~\ref{fig:1} and \ref{fig:2} show the relationship between the value of $\lambda$ and $\theta$ on the expected value of the $Z$ score. The expected value of $Z$ is calculated by randomly drawing 10,000 response outcomes from the distributions given specific parameter values, and then evaluating Equation~(\ref{eq:metric}). More specifically, we fix one of the two variables, and vary the values of the other to assess the impact thereof on $Z$. In Figures~\ref{fig:1} and~\ref{fig:2}, the expected $Z$ score is shown on the $y$-axis, one of the parameters $\theta$ or $\lambda$  is shown on the $x$-axis, and the value of the parameter is depicted in the legend. 

Figure~\ref{fig:1} shows the effect of the response time parameter $\lambda$ on the expected value of $Z$. We first observe the shape of the curve. For small values of $\lambda$, the expected value of $Z$ is close to 0. This is because of response times that are likely to overshoot the maximum permitted response time $d$ (resulting in a direct score of $0$). For higher values of $\lambda$, the expected value of $Z$ grows toward the value of $\theta$. This is due to the fact that higher values of $\lambda$ imply -- on average -- faster response times. In the limit (when $\lambda$ grows very large), the response time does not matter anymore, only the accuracy of the answer (i.e., the upper asymptote is equal to the value of $\theta$). The effect of the parameter $\theta$ is an upward shift of the curve. This is also explained by the previous argument. The logarithmic shape is observed due to the relationship between the threshold for the time and the response time.

Figure~\ref{fig:2} shows the effect of the accuracy parameter $\theta$ on the expected value of $Z$. We see that higher values of $\theta$ lead to a linear increase in the expected value of $Z$. This is to be expected because the definition of $Z$  depends linearly on the accuracy (see Equation~(\ref{eq:metric})). As the value of $\lambda$ increases, the slope of the curve increases up to the point where it models the line $Z = \theta$. This is because $Z$ is penalized less for higher values of $\lambda$. Figures~\ref{fig:1} and \ref{fig:2} verify that the factor of accuracy is more dominant compared to the response speed for the $Z$ score. 

\begin{figure}[!tb]
    \centering
    \includegraphics[width = \linewidth]{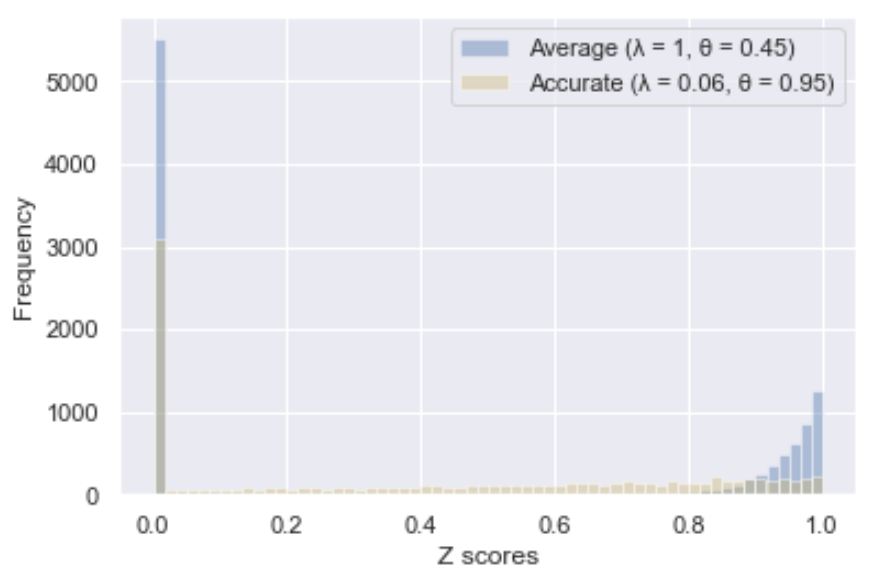}
    \caption{Two opposite mastery profiles that will be attributed with the same expected mastery value of $\mathbbm{E} Z = 0.42$. High accuracy or high speed alone do not necessarily imply a $\mathbbm{E} Z$ higher than 0.5.}
    \label{fig:3}
\end{figure}

Figure~\ref{fig:3} depicts the distribution of $Z$ for two profiles that are both attributed with the same mean value ($\approx 0.42$) for mastery score, but they illustrate two opposed to each other behaviors. The parameters for the profiles are given in the legend. 

The `average' (blue) profile corresponds to a student that provided responses at high speed, and almost half of them are correct. In the figure, that profile has concentrated $Z$ values along the two edges of the $x$-axis. The `accurate' (yellow) profile corresponds to a student that provided almost all answers correctly but took a long time, close to the maximum of the permitted time, to respond. This profile has values of $Z$ that are distributed along with the whole range of values in $x$-axis. We intentionally select these mastery patterns to imply that high accuracy (or high speed) does not necessarily imply a mastery value, which can surpass a mean $Z$ score of 0.5, since the high accuracy is provided at the expense of high response times and the other way around. An expected value of 0.5 would be obtained, for instance, if the `average' (blue) profile would have provided 10\% more correct answers; and for instance if  the `accurate' profile would have responded by 0.05 faster (which is 1/2 of the permitted time). 
A small value for the variance is observed in both cases. For the `average' that is 0.2, and for the `accurate' that is 0.1. 

In reality, there are many things that can influence the accuracy and speed of responses. A loose interpretation for the cause of an `average' profile depicted above is that either the student is using a guessing mechanism (causing an average value for $\theta$) or the skill contains items of control processing \cite{rt2019models}. A similar interpretation can be attributed to the accurate profile, attributing as far more likely the second event. Given similar results regarding the answering speed by other students, we can exclude (or not) the possibility of the second argument. Provided that the possibility for control processing items is excluded, a different interpretation regarding the accurate profile is that it may demonstrate a student that has mastered the skill; yet he acquired it using slow solving strategies \cite{hofman2018fast}. 

\subsection{Convergence rate of the assessment}
\begin{figure*}[!tb]
    \centering
    \includegraphics[width = \linewidth]{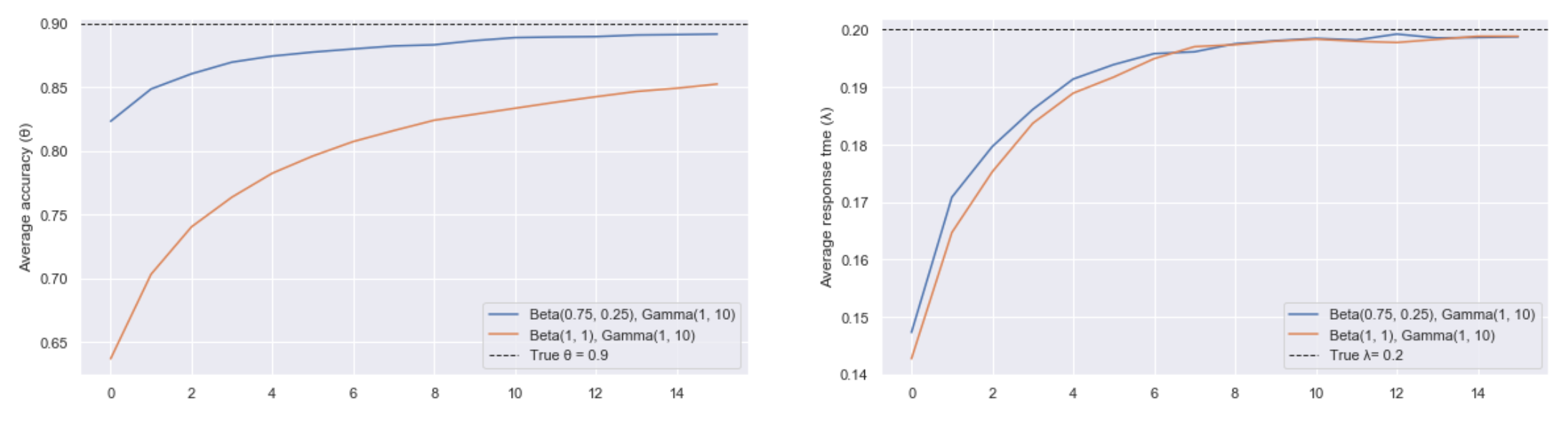}
    \caption{The evolution of the estimated values for $\theta$ and $\lambda$ separately, given an informative (upper line in the left and lower line in the right) and uninformative prior. }
    \label{fig:4}
\end{figure*}

\begin{figure}[!b]
    \centering
    \includegraphics[width = \linewidth]{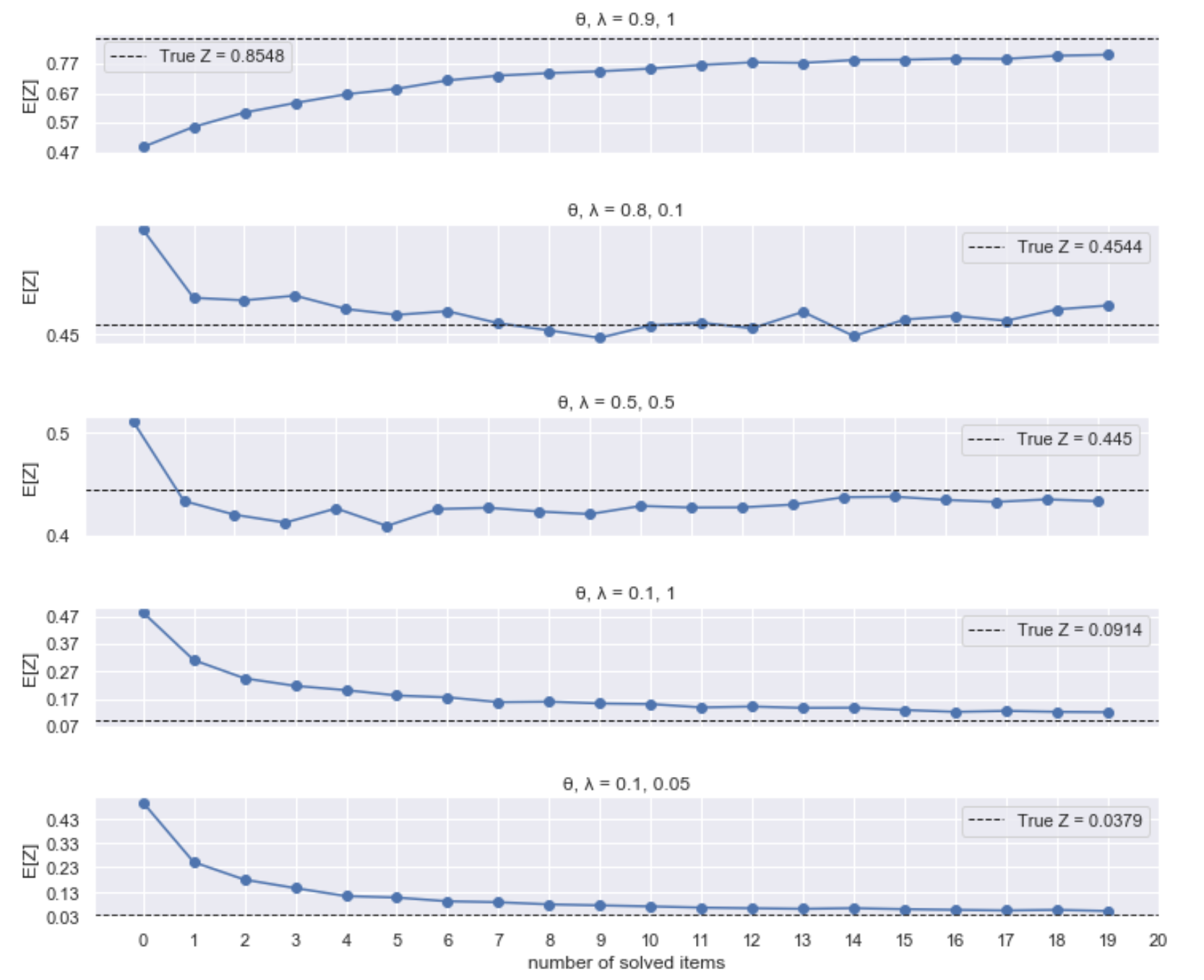}
    \caption{The estimate of $Z$ converges quickly over a sequence of exercises for a set of different profiles with an uninformative prior. The error, the difference between estimated and the true $Z$, is getting smaller after 5 responses.}
    \label{fig:5}
\end{figure}

\begin{figure}[!b]
    \centering
    \includegraphics[width = 9cm]{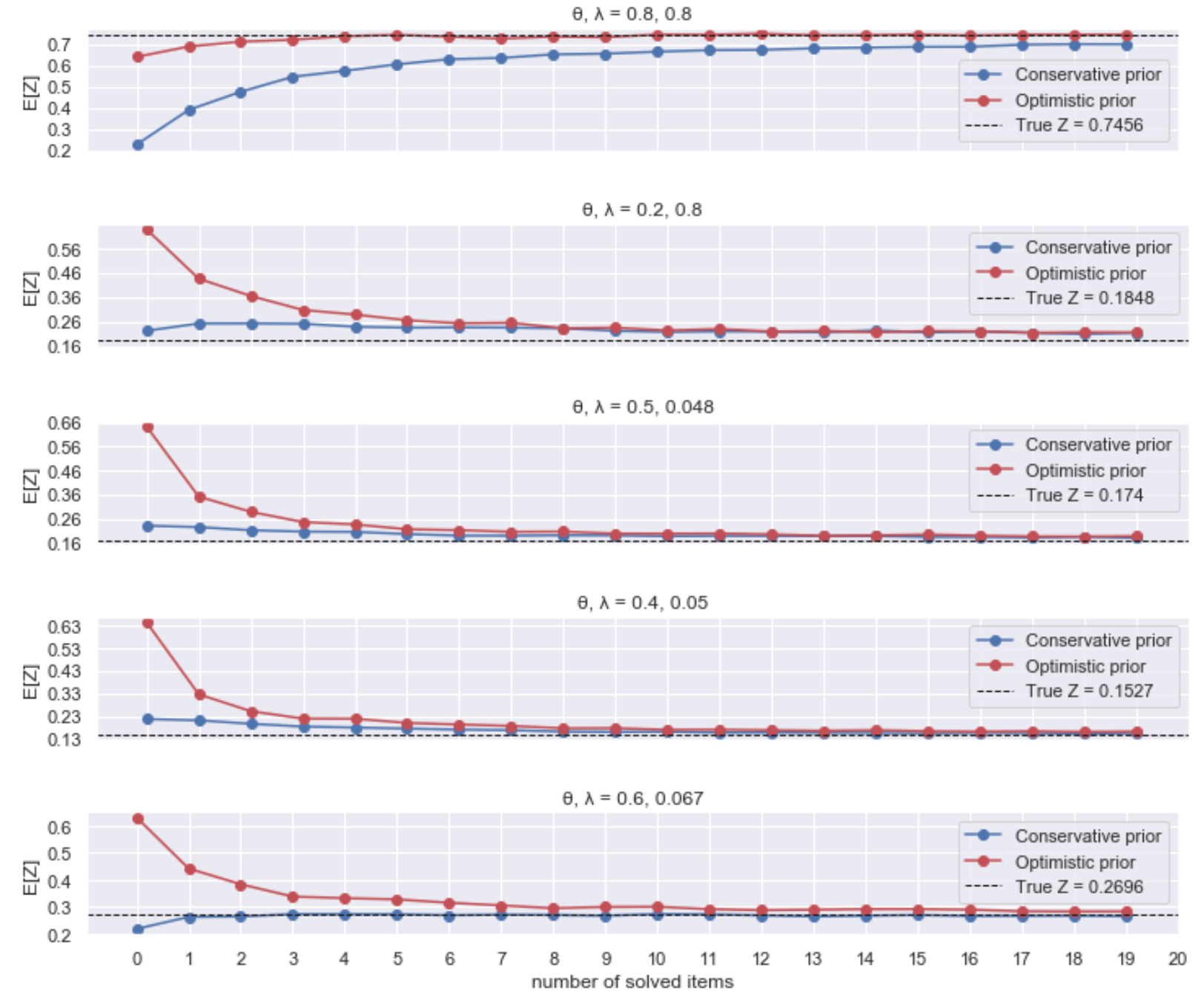}
    \caption{The convergence of $Z$ over a sequence of exercises for an uninformative conservative prior Beta($\alpha=1$, $\beta = 1$), Gamma($n = 1$, $\gamma=10$) and an informative optimistic prior Beta($\alpha=0.75$, $\beta = 0.25$), Gamma($n = 1$, $\gamma=1$) . The informative prior converges sooner only when it aligns with the real mastery level.}
    \label{fig:6}
\end{figure}

In this experiment, we are simulating students who are sequentially providing answers to items. Each time an answer is provided, we examine the change of the estimated parameters $\alpha, \beta, n, \gamma$. As the assessment progresses, we evaluate how many items are needed for the model to learn the true mastery profile (values of $\theta$, $\lambda$, and thereby $Z$). The results of this experiment are as depicted in Figures~\ref{fig:4}, \ref{fig:5}, and~\ref{fig:6}. The dashed black line in Figure~\ref{fig:4} depicts the true parameters $\theta$ and $\lambda$; in Figures~\ref{fig:5} and~\ref{fig:6} it depicts the true $Z$ value. The number of items that the student has answered is always depicted on the $x$-axis.

In Figure~\ref{fig:4}, we show the speed at which the algorithm picks up the value of the true parameters of accuracy $\theta$ and the response time $\lambda$ separately over items. For a mastered skill, represented as a combination of real values of $\theta = 0.9$ and $\lambda= 0.2$, the evolution of the values for the parameters are shown with two lines, each corresponding to a different prior with respect to the accuracy within each plot. The upper line corresponds to the uninformative prior of a Beta($\alpha = 1$, $\beta = 1$) and the lower line to the informative prior of a Beta($\alpha=0.75$, $\beta = 0.25$). One might choose such an informative prior signifying that the estimate for the accuracy is higher than $0.5$, assuming that students are more likely to provide correct answers. We consider that as reasonable given a mastery assessment, a process which requires that students have practiced and learned the skill up to a level before they take the assessment. 

The prior distribution for the response time $\lambda$ parameter is chosen to be an uninformative prior as Gamma($n = 1$, $\gamma=10$) distribution such that the parameter $\gamma$ is equal to half the value of the threshold concerning a single response (i.e., $ n = 1, \gamma = d/2$). An informative prior can lead to faster convergence as it aligns with the true mastery level ($\theta = 0.9$). With an uninformative prior, the model approaches a combination of $\theta= 0.84$ and a $\lambda = 0.2$ within the first 10 exercises. This shows that even with an uninformative prior, the model quickly, with fewer than 10 items, adapts to the true parameters with a low accuracy error (less than 0.1).

Figures~\ref{fig:5} and~\ref{fig:6} show the convergence of $Z$ over items and the sensitivity to the prior for the 5 different student profiles (fixed combinations of $(\theta, \lambda)$). We intentionally choose to depict a range of different student profiles. For the illustration of the experiment, we classify them as follows: (i) fluent $(\theta \ge 0.75, 0.5 \le \lambda \le 1)$ ), (ii) accurate ($\theta \ge 0.75, \lambda < 0.5$), (iii) average ($0.25 < \theta < 0.75$, $\lambda < 0.5)$, (iv) wheel-spinners or unengaged ($\theta \le 0.25$).  The value of $\lambda = 0.5$ is selected by placing the values of the slowest ($d=20$) and fastest ($\lambda=1$) response in the reward function we previously described, and computing the mean distance between the two. Note that, these values were set only for the exposition of showing the model's behavior. In reality, they could take any value and they are set informally by the stakeholders (e.g., tutor) of the algorithm. This allows to explore their true $Z$ values and the convergence of the estimate for $\mathbbm{E}[Z]$. The estimate of the expected value of $Z$ is shown on the $y$-axis.

Figure~\ref{fig:5} corresponds to the case with the uninformative prior where the assessment starts with a prior mastery belief of  $Z_0 = 0.5$. We assume an initial case of a student who is able to solve 3 items at an average response speed (i.e., $10$ seconds) with 95\% probability of correctly answering all of them. That is setting the parameters of the prior to the values of  $\alpha = 0.95, \beta = 0.05, n = 3, \gamma = 30$. 

We observe that the expected value of $Z$ grows closer to the values resembling the true score for all the students as they provide more answers. In general, at the start, the model has a relatively high accuracy error, but as the assessment proceeds, this grows smaller. Over the first five items, the variance is small (i.e., the curve does not highly fluctuate) and hence a smooth convergence is observed. That is especially true for the first (a fluent) and two last plots (wheel-spinners or low-engaged) in contrast to the second (accurate) and third (average) profile. The latter pattern occurs as the two parameters of accuracy and response time are not agreeing in the sense of both taking either high or both taking low values. This results in $Z$ to get values from 0.2 towards 0.5. The variance of the estimate gets smaller as more than 5 items are provided. Over the first to seven items, the accuracy error is small (i.e., the curve reflects the real value with an error less than 0.1) depending on the starting prior and the student profile. Even with an uninformative prior, within the first seven exercises, we argue that the error (less than 0.1) and the variance between the estimate and the real value of $Z$ is relatively small for all profiles. 

In Figure~\ref{fig:6}, we compare an informative and an uninformative prior to the Beta and Gamma distribution. We define an optimistic informative prior that starts with an initial $\mathbbm{E} Z = 0.65$. This profile has a Beta($\alpha=0.75$, $\beta = 0.25$) and Gamma($n = 1$, $\gamma=1$). The uninformative prior starts with a Beta($\alpha=1$, $\beta = 1$) and Gamma($n = 1$, $\gamma=10$). We call the latter a conservative prior as it results in an initial value of $\mathbbm{E} Z = 0.25$. 

We compare the efficiency of convergence between these two different priors for a different new set of student profiles. Specifically, we slightly decrease the accuracy for the `correct' and `incorrect profiles and slightly decrease the response times for all of the five profiles so as to end with one fluent, one unengaged, and three average student profiles.  We intentionally depict these profiles as we want to show four educational cases where the uninformative prior (with regard to the two prior distributions Beta and Gamma) can be closer to the true value of $Z$.

The conservative prior is represented by the first line in the legend (blue) and the optimistic by the second line (red). In the first plot, we have a student who responds at a fast speed and with high accuracy. In that case, the optimistic prior will converge faster because it aligns well with the real mastery level. In contrast to the rest of the profiles, the prior is not necessarily representative of the true mastery profiles, and the uninformative prior converges sooner. Yet, even within the first five items, both of the priors converge at a relatively similar rate as at that point of the assessment both demonstrate an accuracy error less than 0.1 and low variance. In general, starting with the uninformative prior provides a slighter safer choice and treats all student profiles equally. Notwithstanding are the third and last plot (e.g., average profiles). A student demonstrating an average accuracy of 0.5 or 0.6 at a speed slightly less (0.048) or slightly higher (0.067) than the threshold will be attributed in reality with a low mastery score (i.e., $\mathbbm{E} Z=0.17$ or $\mathbbm{E} Z =0.26$). 

Figures~\ref{fig:5} and~\ref{fig:6} show that when the true value of $Z$ lies far from 1 or 0 (i.e., $Z = 0.3$), we observe a higher variance.  This algorithm has fast convergence for all individuals; the ones who have clearly mastered the skill, or have clearly not mastered the skill, and even for wheel-spinning students. On the other hand, the algorithm has slightly slower convergence for those individuals for whom the mastery decision is not as clear, i.e., profiles that demonstrate performance that is highly accurate and really slow, or demonstrate performance that is average both in accuracy and response speed. Regardless of the starting prior and the true student profiles, we argue that within five items, we can recover the true profile given consistent student answers. In addition to that, starting with an uninformative prior for the Beta and the Gamma distribution may be considered a safe choice. That is especially true when there is a broad range of mastery levels among students that are taking the assessment.

\begin{figure}[!t]
    \centering
    \includegraphics[width = 9cm]{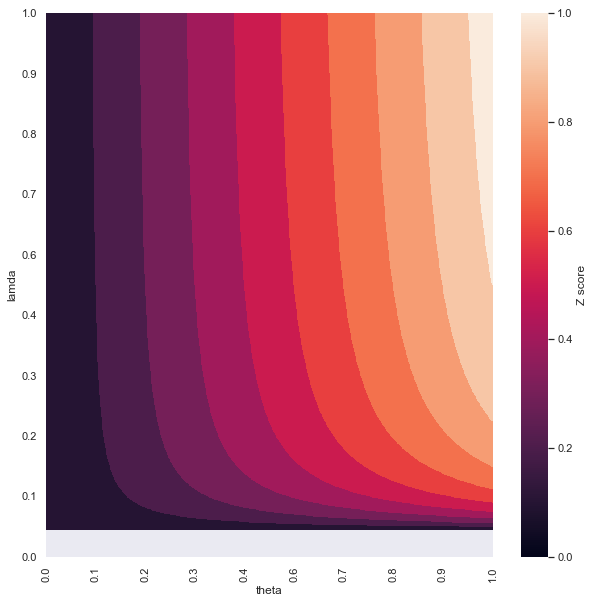}
    \caption{A $Z$ score value is matched to many $\theta$ and $\lambda$  values.}
    \label{fig:7}
\end{figure}

\subsection{Stopping criterion}
In this experiment, we provide argumentation to determine the threshold $\xi$ that will be used as a stopping criterion. Figure~\ref{fig:7} shows isocurves of the $Z$ score (separated in intervals of length $0.1$) for different values of $\theta$ and $\lambda$. Clearly, one cannot simply divide the ranges of the accuracy and responsiveness parameters into intervals to distinguish different student profiles, since the isocurves are not linear. Thus, the $Z$ score cannot be interpreted as the expected portion of correct exercises and the expected response times, and be captured with this classification in a simple table. Instead, rules within $Z$ intervals can be found. A rule could be for instance, a proficient student when he/she has a combination of response times and accuracy that results in a high expected value of $Z$ with adequate certainty. 

Based on the results from the previous experiments and the literature, we have gathered a set of rules for an assessment in Table I. We additionally set the threshold to a value of $\xi = 0.5$ for binary mastery decisions. To consider the students who will not be able to surpass that threshold (i.e., not fluent), similar to previous research, we argue as a valid method to measure whether the $Z$ value gets smaller over the progress of the assessment. The results from the previous experiment indicate that the first three items is a sufficient number of items for that. 

\newcommand{\cmark}{\ding{51}}%
\begin{table}[!tb]
\caption{Mastery profiles}
\begin{center}
\begin{tabular}{cccl}
$Z$ range& Mastery profiles & Decision\\
\hline \hline
0 -- 0.09 &  Wheel-spinning & Remedial strategy\\
0.1 -- 0.25 &  Unengaged & Remedial strategy\\
0.26 -- 0.49 &  Undetermined & Weak Go back\\
0.5 -- 0.74 & Familiar & Weak Pass\\
0.75 -- 0.84 & Proficient & Pass\\
0.85 -- 1 & Mastered &  Strong Pass\\
%\hiderowcolors stop alternating row colors from here onwards\\
\hline  % Please only put a hline at the end of the table
\end{tabular}
\end{center}
\label{table:profiles}
\end{table}

There are several ways to implement that stopping policy, such as expecting a reduction of the $Z$ value within 3 items. Specifically, we pick an optimistic prior that assumes a familiar with the skill student with $\mathbbm{E} Z = 0.5$ (the priors are given by a Beta($\alpha=0.95$, $\beta = 0.05$) and Gamma($n = 3$, $\gamma=30$)). Given a wheel-spinning or an unengaged student, the expected $Z$ value will get reduced after each response. Then it is quite clear that the student will not reach the mastery threshold of $\xi$ and should stop the assessment allocating their learning time differently.  

\begin{figure}[!tb]
    \centering
    \includegraphics[width = 9cm]{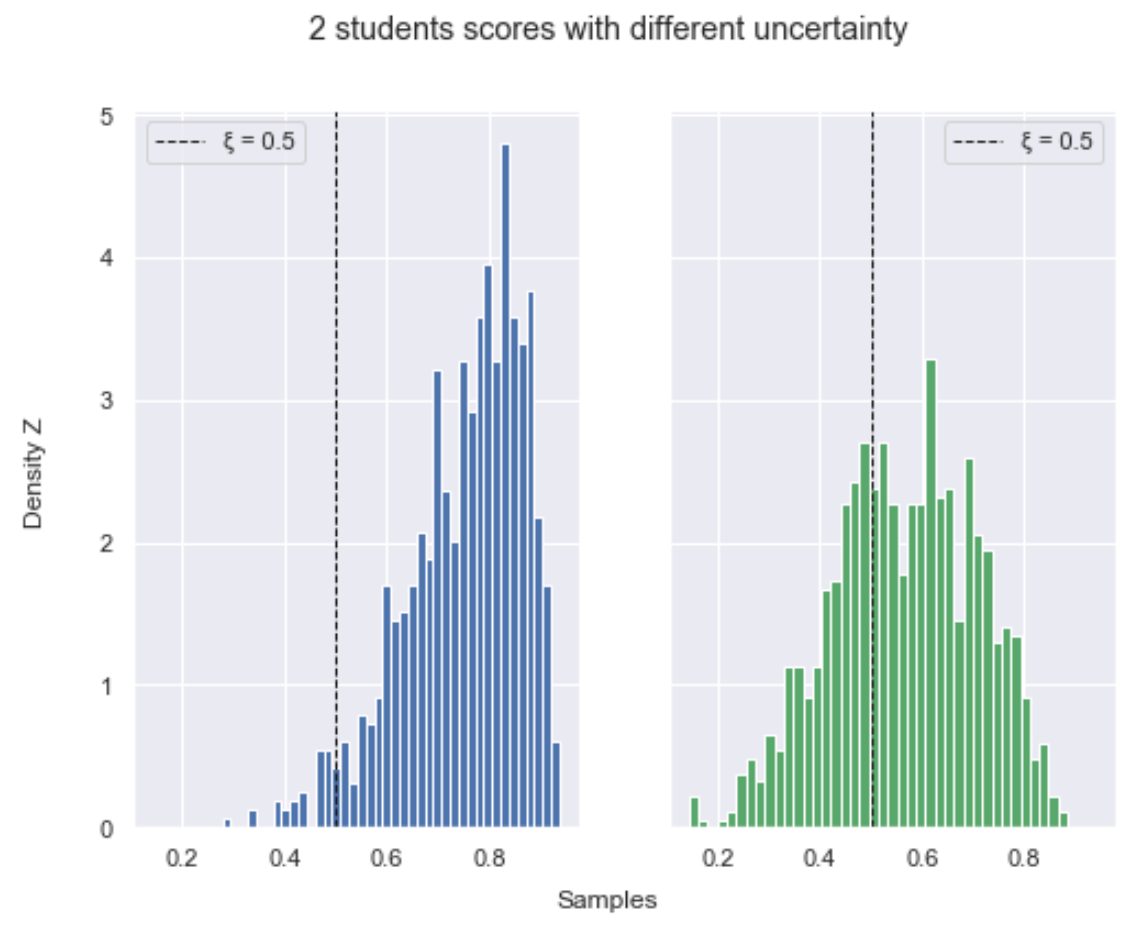}
    \caption{Two students demonstrating similar prior practicing sequences whose expected value for mastery have surpassed the mastery threshold $\xi =  0.5$ (left is attributed with $\mathbbm{E}[Z] = 0.75$ and right with $\mathbbm{E}[Z] = 0.56$ ). Due to the use of a distributional mastery estimate, they will receive a different stopping decision.}
    \label{fig:8}
\end{figure}

With the experiment shown in Figure~\ref{fig:8}, we illustrate the added value of retaining the probability distribution of $Z$, rather than using a point-based mastery estimate. We set the predetermined mastery level to $\xi = 0.5$, and we simulate two students with similar prior practicing sequences. They have both solved 10 questions at an average speed of 0.5 seconds, and they only differ in that one has answered two more questions incorrectly compared to the other. More formally, the state for the one, say $s_1$, is $s_1 = (\alpha = 8, \beta = 2, n = 10, \gamma = 12)$ and for the other, say $s_2$, is $s_2 = (\alpha = 6, \beta = 4, n = 10, \gamma = 12)$. Then the expected score can be calculated for the former as $Z(s_1) = 0.7$ and for the latter $Z(s_2) = 0.6$. Since their point-based estimates for their scores at this round of the assessment has just surpassed the chosen $\xi$, the optimal decision for both should be to stop the assessment. However, that does not hold when the whole distribution is used as a score due to the high uncertainty surrounding the estimated value of $Z(s_2)$. 

To illustrate the above, we need to evaluate Equation~(\ref{eq:vi}) for the different states $s_1$ and $s_2$. Suppose that the structure of the value function $V(\cdot)$ is linear in its parameters, such as the following $V(\alpha, \beta, n, \gamma) = k_1 \alpha  + k_2 \beta + k_3 n + k_4 \gamma$, with $(k_1, k_2, k_3, k_4) = (1, 1, -0.71, -0.37)$. By evaluating the integral in Equation~(\ref{eq:vi}) with respect to time $t$, we end up with two expressions, one for each student. To decide if a student is attributed with mastery, we use the inequality $V(s) \ge \xi/(1-\eta)$. Therefore, we have a system of two inequalities, each for a student. We additionally set the $V(s_1) > V(s_2)$ and the value for $\eta = 0.8$. We end up with a $V(s_2) = 2.35$ and a $V(s_1) = 2.5$, which implies that the optimal decision for $s_2$ would be to continue the assessment since our algorithm selects the minimum, while the action in state $s_1$ would be to stop. 

It is important to note that in practice, there is the possibility that a student with a true mastery highly close to the threshold may not take the assessment at all as their profile is estimated as `mastered'. That does not sound interesting for the current study, but it would be beneficial for an assessment that involves multiple skills as the student will allocate their assessment time to only the skills that have not been mastered yet.

Our model can be adapted to different assessment purposes. A more certain or a quicker assessment is tuned by setting the value of the hyperparameter $\eta$. A more lenient or stricter assessment is set by the value of $\xi$ while the weight of the desired response time for the skill is set by the threshold $d$. Last but not least, given prior student data, more informative starting priors for individual students can be set such that their progress can be measured more accurately over the time they are taking the assessment.

\section{Conclusion} \label{sec:conclusion}
To conclude, we introduce a model that estimates whether an individual has attained mastery of a skill based on a sequence of responses that leverage the information of speed and accuracy of the previously given answers. The properties of the model are established using numerical experiments. We found that we can recover the mastery profile with low variance within approximately 5 responses. That holds in the case of an uninformative, optimistic or conservative prior belief about student mastery indicating that our model is not sensitive to the prior. We additionally found three results that are aligned with previous research. Firstly, the variance of the length of the assessment is increased for profiles for which the mastery level is not clear. These constitute students whose responses are provided at really low speed and very high accuracy; or at an average speed and accuracy. The number of 5 consecutive responses as sufficient assessment length is a second finding observed in literature. At the same time, shorter assessments are materialized for clearly mastered or wheel-spinning students. Lastly, we explored rules on mastery score ranges, which are combinations of values for response time and accuracy -- not fixed values for these two -- that point to a stopping decision. To resolve for a suboptimal stopping decision, we found that expecting increase or decrease of the starting $Z$ within the first three responses is sufficient. There are four advantages of the proposed model: (i) it takes into account the certainty of the estimate, which can guide to more careful decisions, (ii) it is transparent, allowing us to explain the results. It justifies the evaluation of mastery, (iii) the inclusion of informative priors can be done naturally and easily, offering the potential for more efficient assessments, (iv) it is flexible, allowing for extensions and adaptations to different assessment purposes. We plan to extend the model for a multi-skill assessment in the future.

\subsection{Future directions}
There are plenty of future directions for the model and the research community. The most useful enhancement of the model is to include more than one skill set (or item difficulties). An assessment often tests the attained mastery on a composition of learning goals, often termed as multi-stage testing~\cite{mastery1999bayesian}. In that case, the decision to stop the assessment of the current skill and continue to the next is more informative and realistic compared to assuming an individual skill. 

Regarding the efficiency of the assessment, a representative data sample of student responses can be used to set the hyperparameters of the model, such as the prior beliefs for each student mastery level to be close to the truth, leading to a shorter assessment. Similarly, the value of $d$ should (and can) be determined using empirical data, personal goals of students, or set by teachers. In addition to that, setting a maximum number of items stopping rule would lead to a more efficient assessment for students whose real mastery is undetermined or lower than the targeted one, leading to an extra hyperparameter of the model. In~\cite{gong2015towards}, they found that a low score within the first 10-15 exercises provides a hint that the student should spend his time learning the previous skill or the current skill rather than assessing his/her mastery level. However, they defined the decision of mastery as three consecutive correct number of answers in a row, and they use a different framework. Research on stopping criteria~\cite{cat2019comparestop} also agrees that including such a maximum number of items hyperparameter works better in practice.

The information provided by the assessment can be more detailed with a variable assessment length that matches quantiles of $Z$, providing more than just two decisions (e.g., too slow on `these' items). Combining that with the above-proposed enhancements implies that rather than immediately stopping the assessment (once $(1 - \eta) V(\cdot) \ge \xi$), we could potentially increase the certainty by letting the student provide a few more responses to categorize him/her within a more certain profile. This could lead to multiple profiles of mastery. For instance, we could detect a wheel-spinning type of student to stop with fewer exercises. We illustrate such a scenario in Table~\ref{table:profiles}. Moving forth or back should be taken carefully, since such a student may be stuck in a loop of going either `back-forth' or `back-back-back'.

Lastly, it is well known that a dynamic model should be evaluated with online empirical data, i.e., interventions in an online learning platform. A weakness of every model that is based on simulations is that these are a simplification of reality. They are advantageous in that they provide access to the ground truth and a safety zone for understanding the workings of the model. Simulation studies are necessary to determine the appropriate assessment properties (e.g., test length and termination criteria) required to ensure the degree of precision and test efficiency required in operational testing scenarios \cite{catnewstop}. The results can be seen as an opportunity to learn how to choose the stopping criterion (i.e., $\xi$) to fit the underlying application of the model, either that is practicing (with slow growth) or assessing. In the real-data case, testing for data patterns that lead the model to wrong stopping decisions is always interesting since this will give information on updating the model by incorporating the feedback provided by the stakeholders of the model.  

In general, the choice of the stopping criterion in adaptive assessment and practicing is still an open issue~\cite{kaser2016stop, sapountzi2019dynamic, bktthreshold, masteryrupp, pelanek2017experimental, catnewstop}. Although a threshold that adapts to certain mastery levels would offer flexibility, it is not easy to be made. Fundamentally, that choice is a matter of definition, which means that students have attained the mastery level that is required for the targeted skill. The adaptation to the mastery level (e.g., weaknesses/strengths of a student) in an assessment session can be made on a fine-grained scale where everything is quite well-defined, and there is a specific target skill. This contradicts with the ability of the model to generalize and be used in a different context. 

\bibliographystyle{elsarticle-num} 
\bibliography{paper}

\end{document}